%% file: main.tex
\documentclass[a4paper,psamsfonts,12pt]{amsart}
\usepackage{amsmath,amssymb}
\usepackage[totalwidth=480,totalheight=680]{geometry}
\input epsf.tex
\input MyDefAll.tex

\def\sbullet{\,\,{\boldsymbol \cdot}\,\,}
\title[Abresch-Langer Conjecture]{
On the Saddle Point Property of Abresch-Langer Curves under
the Curve Shortening Flow}
\author[Thomas Au]{Thomas Kwok-keung Au}
\address{Dept of Mathematics, The Chinese University of Hong Kong, Shatin, Hong Kong.}
\email{thomasau@cuhk.edu.hk}

\newtheorem{theorem}{Theorem}
\newtheorem{proposition}[theorem]{Proposition}
\newtheorem{lemma}[theorem]{Lemma}
\newtheorem{corollary}[theorem]{Corollary}

\newtheorem*{ALtheorem}{Abresch-Langer Theorem}
\newtheorem*{ALconjecture}{Conjecture}
\newtheorem*{maintheorem}{Main Theorem}

\def\bn{{\boldsymbol \nu}}
\def\hmn{{h_{m,n}}}
\def\kmn{{\kappa_{m,n}}}

\def\hff{h_{\theta\theta}}
\def\hepsilon{{h_\varepsilon}}
\def\kepsilon{{\kappa_\varepsilon}}
\def\tkappa{{\tilde\kappa}}
\def\tildeh{{\tilde h}}
\def\calE{\mathcal E}
\def\calF{\mathcal F}

\begin{document}
\begin{abstract}
In the study of the curve shortening flow on general closed curves,
Abresch and Langer posed a conjecture that the homothetic curves can be
regarded as saddle points between multi-folded circles and some singular
curves.  In other words, these homothetic curves are the watershed
between curves with a nonsingular future and those with singular future
along the flow.  In this article, we provide an affirmitive proof to
this conjecture.
\end{abstract}

\maketitle
\thispagestyle{empty}
\setlength\baselineskip{28pt}
\setlength\parskip{2ex}
\setlength\parindent{0em}
\pagestyle{headings}

\input intro.tex

\input prelim.tex

\input study1.tex

\input study2.tex

\input bib.tex

\end{document}

%

%% file: MyDefAll.tex
\def\ifnextchar#1#2#3{\let\tmpnce=#1%
    \def\tmpnca{#2}\def\tmpncb{#3}\futurelet\tmpncc\ifnch}%
      \def\ifnch{\ifx\tmpncc\tmpnce\let\tmpncd=\tmpnca%
	      \else\let\tmpncd=\tmpncb\fi\tmpncd}
\def\diff#1#2{\dfrac{\d #1}{\d #2}}
   \def\pdiff#1#2{\dfrac{\partial#1}{\partial#2}}

\def\d{\operatorname{d}\!}

\def\meanint{\hskip{.75em}\bar{}\hskip{-.75em}\int}
\def\modulus#1{\left|#1\right|}
\def\monthname{\ifcase\month\or Jan\or Feb\or March\or Apr\or %
    May\or June\or July\or Aug\or Sept\or Oct\or Nov\or Dec\fi}
\def\norm#1{\left\|#1\right\|}
      \def\normp#1_#2{\norm{#1}_{#2}}

\def\homeo#1#2#3
	{#1~\colon~#2~\stackrel\approxeq\to~#3}

\def\O#1{\text{O}(#1)}

\def\meanint{\hspace*{.75em}\bar{}\hspace*{-.75em}\int}

%% file: intro.tex
Let $\gamma_0$ be a given immersed closed plane curve.  We consider the
initial value problem of the curve shortening flow
\begin{equation}\label{eqn-flow}
\left\{
\begin{aligned}
\gamma_t(p,t) &= -\kappa(p,t)\bn(p,t) \\
\gamma(p,0) &= \gamma_0(p)
\end{aligned}
\right.
\end{equation}
where $\gamma(p,t)$ is a family of curves with curvature $\kappa$ and
unit (outward) normal $\bn(p,t)$.  The curve shortening flow for an embedded
closed initial curve was completely characterized by the Grayson
convexity theorem, \cite{Grayson}, and the Gage-Hamilton theorem,
\cite{GH}.  The first asserts that the flow drives
any such $\gamma_0$ to a convex curve while the second says that a
convex curve shrinks to a round point.  However, when $\gamma_0$
has self-intersections, it is easy to see that singularities may
arise in the process.  A typical example is the flow of the cardoid.  The
little loop of the cardoid contracts and develops a cusp when the large
loop still exists.  It is thus necessary to classify the singularity of the
flow.  A natural way of classification arises from the blow-up rate
of the curvature into type~I and type~II singularities (Altschuler
\cite{Al} and a series of important studies of Angenent, for examples, 
\cite{An3} and \cite{AV}).  In our situation, all the singularities
are of type~I, namely,
$ \modulus{\kappa}_{\max}(t) \sqrt{t_\infty-t} \leq C, $
where $t_\infty$ is the blow-up time, the singularity looks
asymptotically like a contracting self-similar solution.

A (contracting) self-similar solution of~(\ref{eqn-flow}) is a flow
in which the shapes of the curves change homothetically and
continuously to a point in finite time.  A curve in the flow is called a
contracting self-similar curve.  Obviously, the circle is a contracting
self-similar curve.  It turns out that other such curves must have
self-intersections. In fact, all closed contracting self-similar
solutions had been completely classified by Abresch and Langer,
\cite{AL}.  For our future discussion in this article, it is
convenient to express their result in terms of the support function
of the curve,
since all these solutions are locally convex.

For each locally convex curve, $\gamma(p)$, the support function is
a function
$$
h(\theta) = \langle \gamma(p), \bn(p) \rangle,
$$
where $\theta$ is the angle of the outward
normal and $\bn=(\cos\theta,\sin\theta)$.
The position of the curve $\gamma$ in terms of $h$ is given by
$$
\gamma(p) = \left(
h(\theta)\cos\theta - h_\theta(\theta)\sin\theta,
h(\theta)\sin\theta + h_\theta(\theta)\cos\theta \right). 
$$
This formulation is common in the discussions of convex plane
curves.  For reference, one may see \cite{C}.  Using the formula
$$
\kappa(p) = \frac{1}{h(\theta)+\hff(\theta)},
$$
equation~(\ref{eqn-flow}) can be rephrased as
$$
h_t(\theta,t) = \frac{-1}{h(\theta,t)+\hff(\theta,t)}
$$
and the support function $h(\theta)$ of a contracting self-similar
solution satisfies
\begin{equation}\label{eqn-hmn0}
\hff + h = \frac{C}{h}, \qquad\quad h > 0, \quad C > 0,
\end{equation}
where $C$ is a positive constant.  Without loss of generality, from
now on, we take $C=1$.  Moreover, equation~(\ref{eqn-hmn0}) has a
finite integral
$$
h_\theta^2 + h^2 = 2\log h + C_0,
$$
It is easy to see that its solutions are positive and periodic.  
Let us denote by $h(\theta; \alpha)$ the
solution to~(\ref{eqn-hmn0}) satisfying the initial conditions
$$
h(0;\alpha) = \alpha > 1 \qquad \text{and} \qquad h_\theta(0;\alpha) = 0.
$$
Then the main result in~\cite{AL} can be stated as follows:
\begin{ALtheorem}
First, the circle is the only embedded contracting self-similar curve.
Second, as $\alpha$ increases from~1 to~$\infty$, the period of 
$h(\sbullet;\alpha)$ decreases strictly from
$2\pi/\sqrt{2}$ to~$\pi$.
\end{ALtheorem}
Whenever the curve determined by the support function
$h(\sbullet;\alpha)$ is closed, its
period must be a multiple of~$2\pi$.  In other words, for any pair
of relatively prime positive integers $m, n$ satisfying $\frac{1}{2} <
\frac{m}{n} < \frac{1}{\sqrt{2}}$, there corresponds a unique
contracting self-similar curve with $n$~leaves in $m$~rotations.
In this article, we follow the notation of Abresch-Langer to
denote such a homothetic curve by $\gamma_{m,n}$ and, likewise,
its support function by $\hmn$.
\begin{center}
\mbox{\epsfysize=49mm \epsfbox{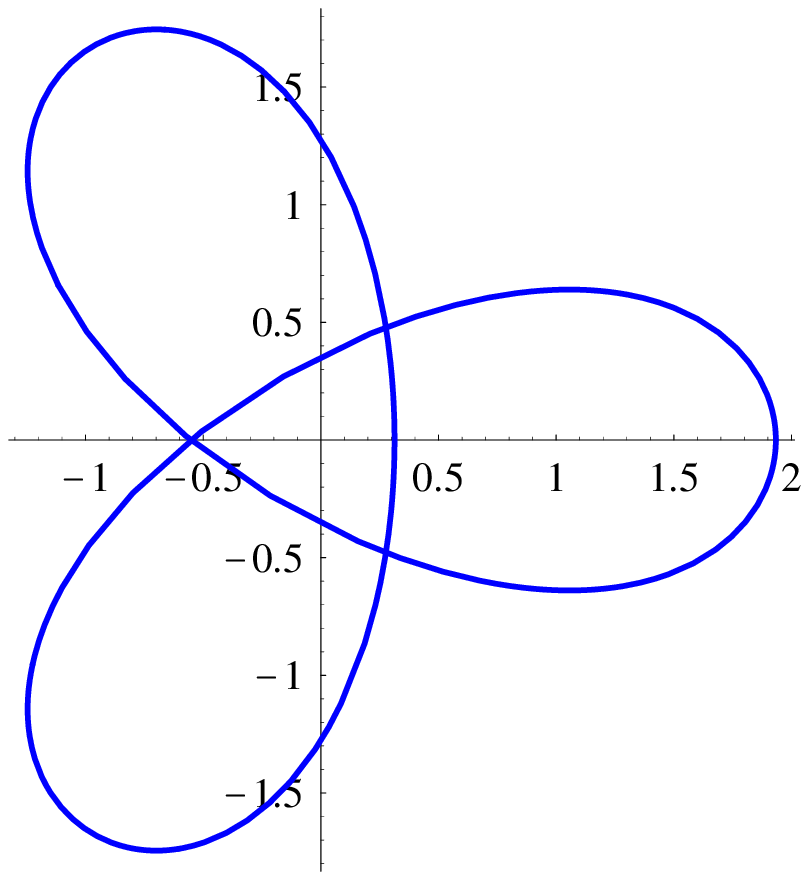}}\hfil
\mbox{\epsfysize=49mm \epsfbox{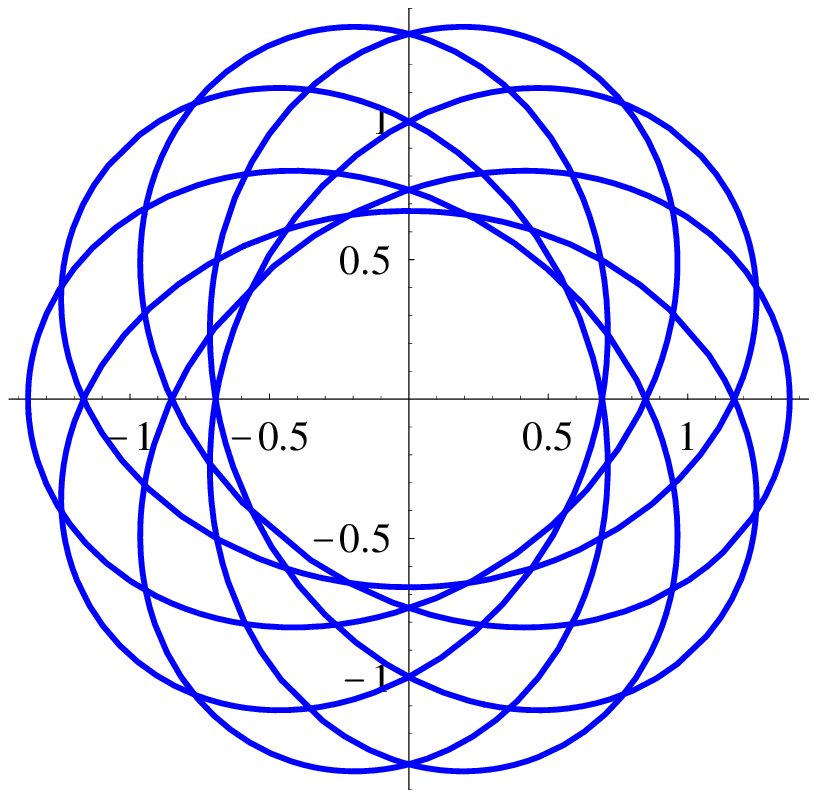}}\hfil
\mbox{\epsfysize=49mm \epsfbox{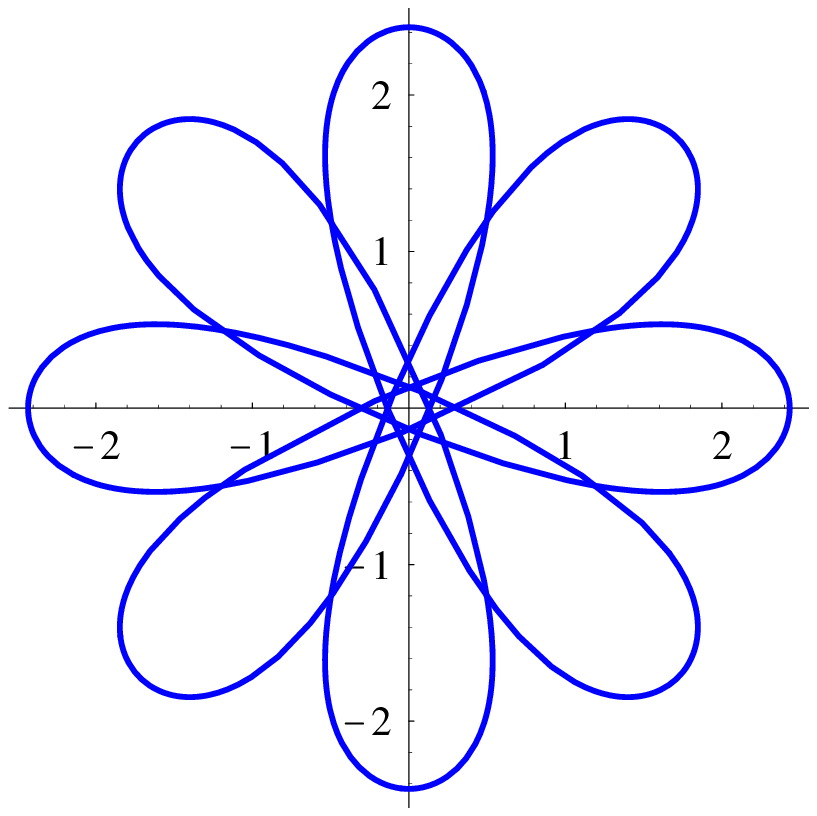}} \\
\nobreak
\mbox{\footnotesize
Figure~1: pictures of $\gamma_{2,3}$, $\gamma_{7,10}$, and
$\gamma_{5,8}$.
}
\end{center}
Geometrically, the value $\alpha$ is the distance from the origin
to the tip of a leave; and analytically, it is the maximum of $h$.

The linear stability properties of these contracting self-similar
solutions are studied
in \cite{AL} and Epstein-Weinstein, \cite{EW}.
As related to the nonlinear stability of these curves, Abresch
and Langer explain that they serve in a certain sense as
``saddle points'' between nonsingular and singular curves.  More
precisely, their conjecture can be stated as follows.
\begin{ALconjecture}
Consider the equation~{\em(\ref{eqn-flow})} with
initial data $\gamma_0 =
\gamma_{m,n}+\varepsilon\bn$, with $\modulus{\varepsilon}$ small.
\begin{enumerate}
\item[(a)]
When $\varepsilon>0$, the trajectory through $\gamma_0$ is
asymptotic to an $m$-fold circle; and
\item[(b)]
when $\varepsilon<0$, the trajectory through $\gamma_0$ is
asymptotic to a singular curve $\Gamma_{m,n}$ with $n$~cusp
points.
\end{enumerate}
\end{ALconjecture}
They also expect the evolution process after rescaling behaves
as shown in the following pictures.
\begin{center}
\begin{tabular}{ccc}
\mbox{\footnotesize
Figure 2a: rescaled evolution for $\varepsilon>0$.}
&\hbox{\hspace*{1cm}}&
\mbox{\footnotesize
Figure 2b: rescaled evolution for $\varepsilon<0$.} \\
\mbox{\epsfysize=49mm \epsfbox{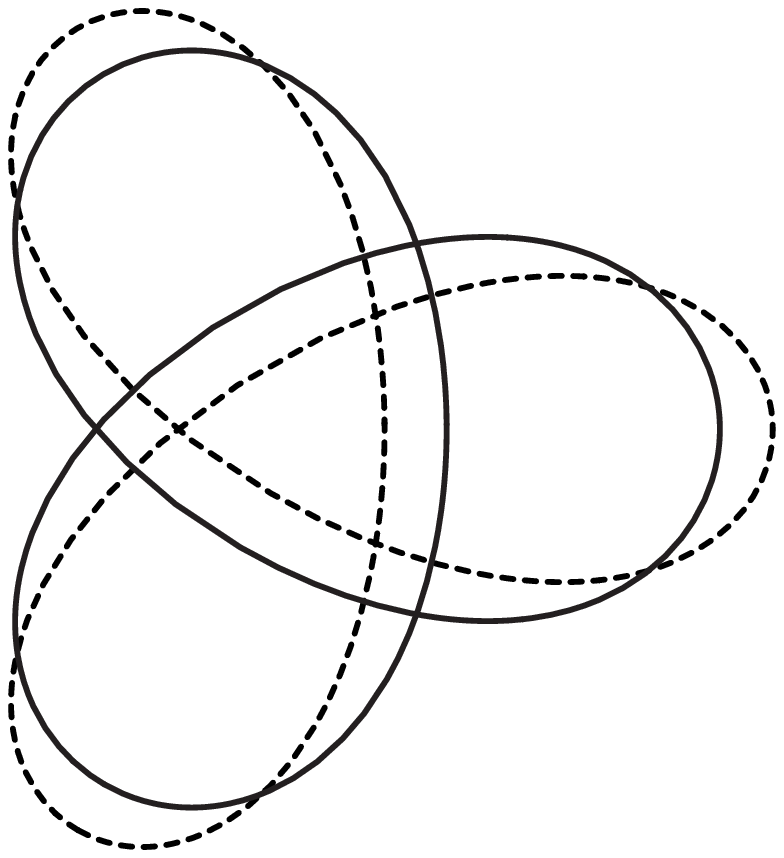}} 
&\hbox{\hspace*{1cm}}&
\mbox{\epsfysize=49mm \epsfbox{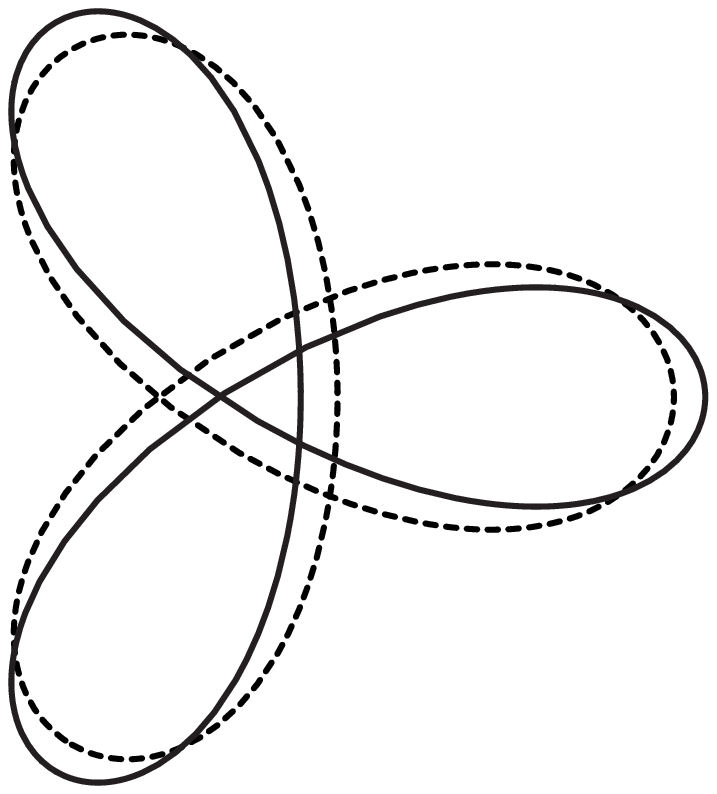}} \\
\mbox{\epsfysize=48mm \epsfbox{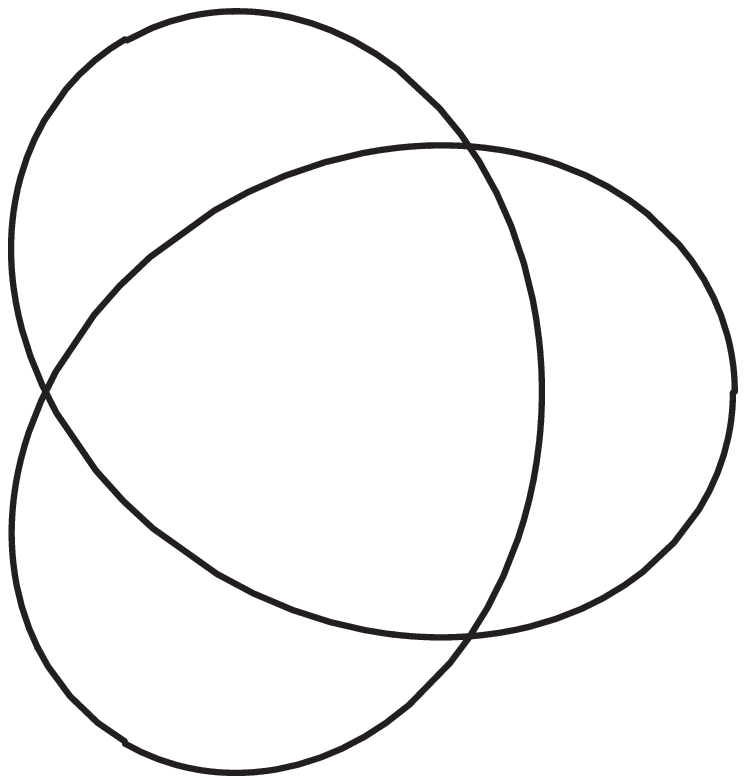}} 
&\hbox{\hspace*{1cm}}&
\mbox{\epsfysize=48mm \epsfbox{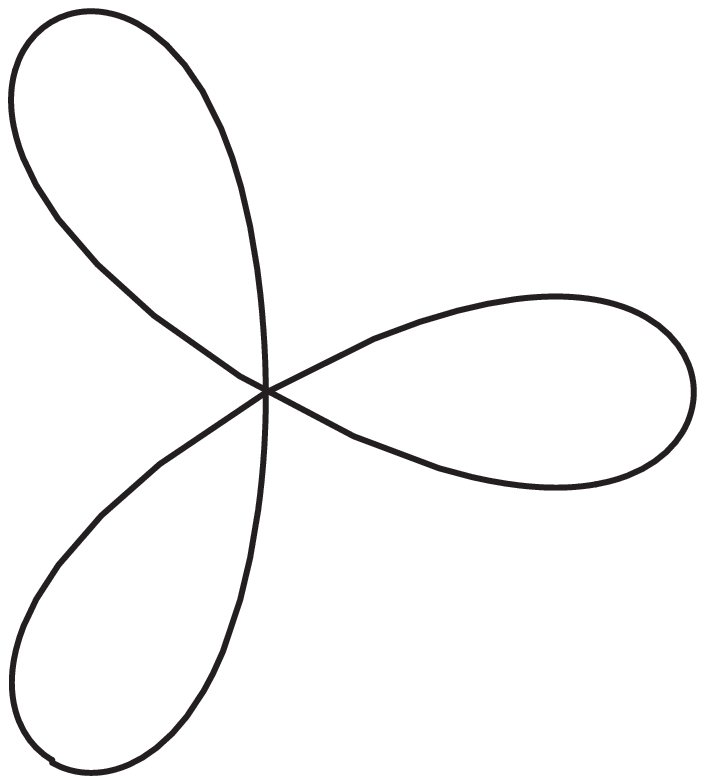}} \\
\mbox{\epsfysize=48mm \epsfbox{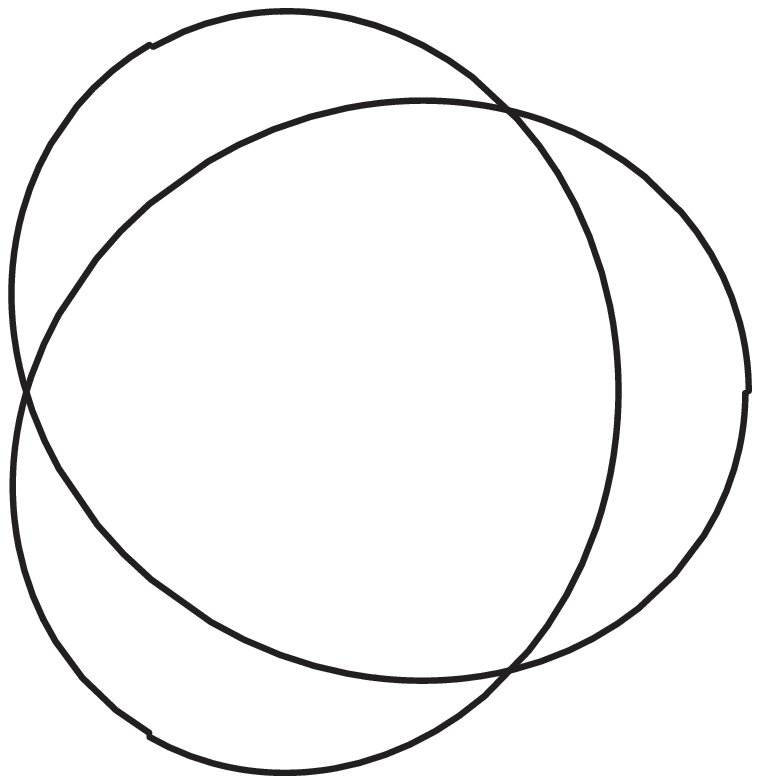}} 
&\hbox{\hspace*{1cm}}&
\mbox{\epsfysize=48mm \epsfbox{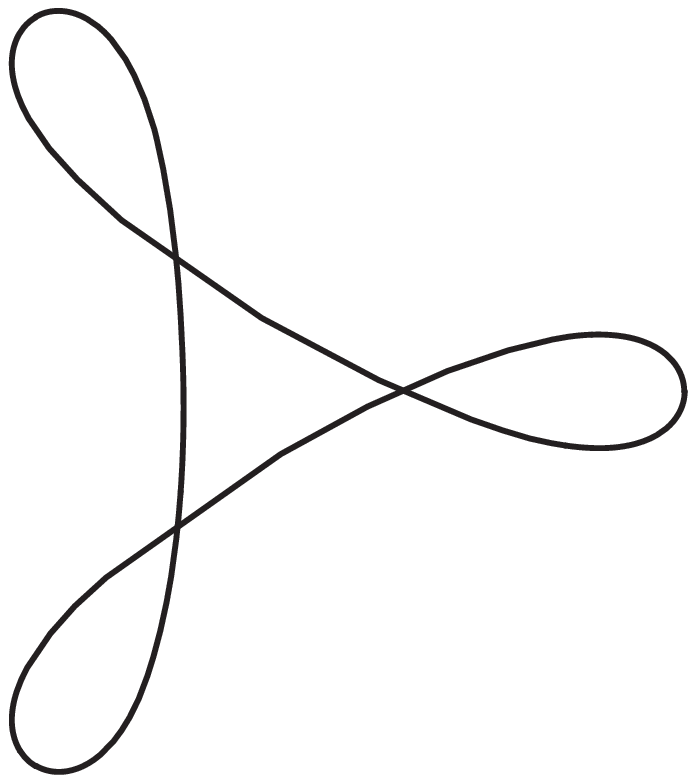}} \\
\mbox{\epsfysize=49mm \epsfbox{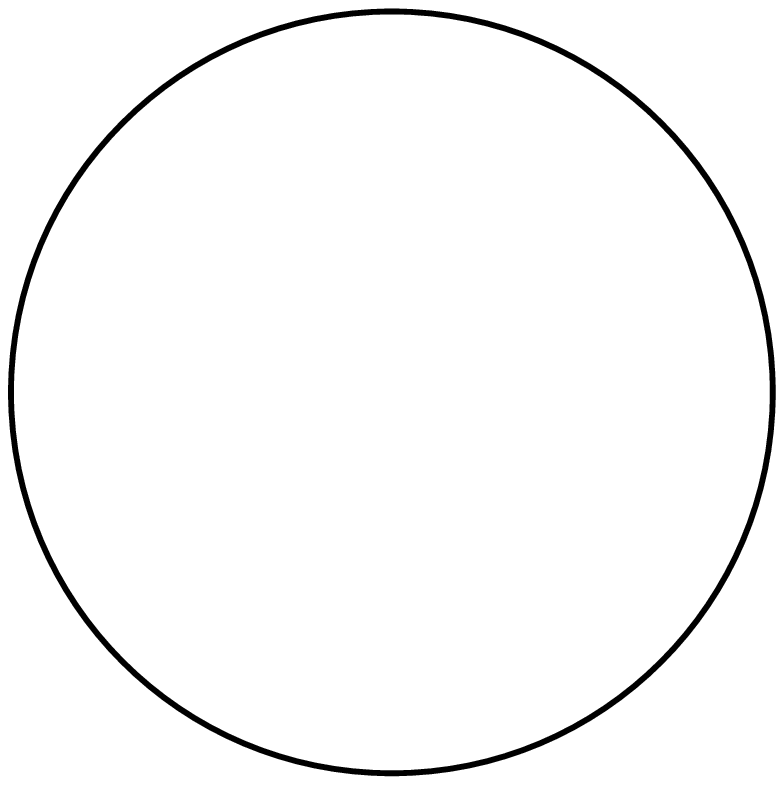}} 
&\hbox{\hspace*{1cm}}&
\mbox{\epsfysize=49mm \epsfbox{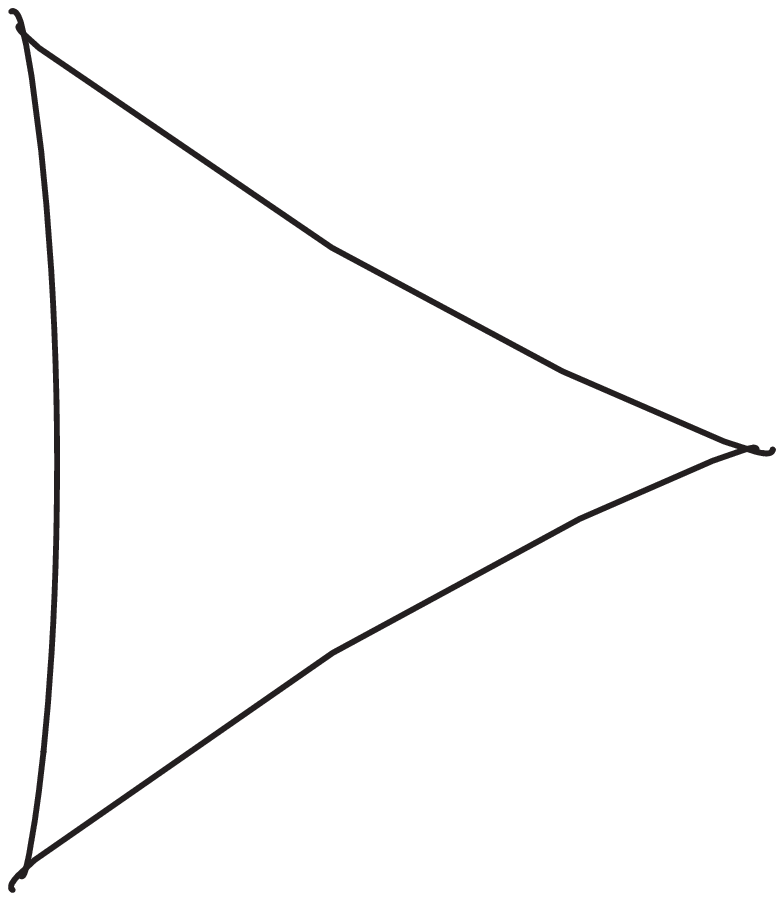}} \\
\multicolumn{3}{c}{\parbox{140mm}{\footnotesize
These pictures are numerically produced with a Mathematica
algorithm which is mostly according to the flow with
minor modification to overcome the slow convergence of the
program.  All the curves have been rescaled.
The leftmost picture shows  homothetic curve
$\gamma_{2,3}$ in dash and its perturbation $\gamma_{2,3}+\varepsilon\bn$
in solid.}}
\end{tabular}
\end{center}

We shall prove the conjecture in this article.   Furthermore, we
also show analytically how the flow develops singularities in
case~(b) and the shape and structure of the curve is the same as
the qualitative description in \cite{AL}.
\begin{maintheorem}
The Abresch-Langer conjecture is true.  
\begin{enumerate}
\item[(a)]
In the case that
$\varepsilon>0$, solution to~{\em(\ref{eqn-flow})} exists for all
time $t\in[0,\infty)$ and it tends to an $m$-circle as
$t\to\infty$.  
\item[(b)]
In the case that
$\varepsilon<0$, $n$~singularity points are formed when the area of a
leaf becomes zero and the curvature blows up.
\end{enumerate}
\end{maintheorem}
Our discussion is organized in the following way.  In \S 1,
we formulate the problem in
a suitable and convenient way.  Several equivalent versions of the
problem will be present; each of which will be used in later
discussion.  The formulation of the problem in terms of
support function has been given in \cite{C}.  In fact,
such formulation provides a clean and elegant view of the
topic, which is essential to our solution of the problem.
In \S 2, we first lay the foundation for applying maximum
principle to the solution.  It follows naturally with the
derivative estimates and the convergence of the flow to
self-similar solution in the case that $\varepsilon>0$.  Then the
instability of Abresch-Langer curves is shown.
In \S 3, we deal with the case that $\varepsilon<0$.  More
specifically, we first show that the solution always exists before the
second stage in Figure~2b.  It further exists as long as the leaves do
not vanish.  Since it is known that these leaves eventually 
eventually shrink to a point,  at this moment, the
singularities occur.

The discussions from our colleagues have been very encouraging and
insightful in our study.  We would like to specially mention
K.~S.~Chou and Tom Wan.

%% file: prelim.tex
\setcounter{section}{0}
\section{Formulation and Normalization}
First of all, let us
reformulate the problem in terms of support functions of the curves.
We would like to study the initial value problem of solving for
$h(\theta,t)$,
$$
\left\{
\begin{aligned}
h_t(\theta,t) &= \frac{-1}{h(\theta,t)+\hff(\theta,t)}, \\
h(\theta,0) &= \hmn(\theta) + \varepsilon.
\end{aligned}
\right.
$$
In order to compare the problem with the situation at
$\varepsilon=0$, we replace $\hmn+\varepsilon$ with
$$
\hepsilon(\theta) = \left(1 + \frac{L}{m\pi}\varepsilon +
\varepsilon^2 \right)^{-1/2} \left(\hmn(\theta)+\varepsilon\right),
$$
where $L$ is the arc length of $\gamma_{m,n}$.  Next, the algebraic area
of the curve $\gamma$ can be determined by the support function $h$,
namely,
$$
A(h) = \frac{1}{2}\int_0^{2m\pi} h(h+\hff) \d\theta.
$$
Then, $A(\hmn) = A(\hepsilon) = m\pi$ with the above choice of
initial support function.  Let $h(\theta,t)$ be the
periodic solution to the evolution equation for $\theta \in
[0,2m\pi]$ and $t\geq 0$,
\begin{equation}\label{eqn-flow0}
\left\{\begin{aligned}
h_t(\theta,t) &= \frac{-1}{h(\theta,t) + \hff(\theta,t)}, \\
h(\theta,0) &= \hepsilon(\theta). \\
\end{aligned}\right.
\end{equation}
We then have
\begin{align*}
\pdiff{A(h(\sbullet,t))}{t} &= \frac{1}{2}\int_0^{2m\pi} h_t(h+\hff) +
h(h_t+h_{t\theta\theta}) \d\theta \\
&= \int_0^{2m\pi} h_t(h+\hff) \qquad \text{after twice integrating
by parts,} \\
&= \int_0^{2m\pi} (-1) = -2m\pi.
\end{align*}
In other words, $A(h) = A_0 - 2m\pi t = m\pi(1-2t)$.  We further
normalize the flow by letting
$$
\tilde{h}(\theta,t) = \frac{1}{\sqrt{1 - 2t}}h(\theta,t)
$$
and changing of variables from $t$ to $\tau$ such that
$\diff{\tau}{t} =1/(1-2t)$.
We have the following initial value problem for $\theta \in
[0,2m\pi]$ and $\tau\geq 0$,
\begin{equation}\label{eqn-flow1}
\left\{\begin{aligned}
{\tilde h}_\tau(\theta,\tau) &= -{\tilde\kappa}(\theta,\tau) +
{\tilde h}(\theta,\tau), \\
{\tilde h}(\theta,0) &= \hepsilon(\theta). \\
\end{aligned}\right.
\end{equation}
As long as $A(h) > 0$, (\ref{eqn-flow1}) is equivalent to
(\ref{eqn-flow0}).  This is true in particular if $h>0$.  
Under this formulation, $\hmn$ is a stationary
solution to (\ref{eqn-flow1}).

Let $\kappa$ and ${\tilde\kappa}$ denote the curvatures of the curves
with support functions $h$ and ${\tilde h}$ respectively.  Then they
satisfy the following equations, which correspond to
equations~(\ref{eqn-flow0}) and~(\ref{eqn-flow1}) respectively,
\begin{align}
\kappa_t &= \kappa^2 ( \kappa_{\theta\theta} + \kappa ),
\label{eqn-flowk}\\
{\tilde\kappa}_\tau &= {\tilde\kappa}^2 \left({\tilde\kappa} +
{\tilde\kappa}_{\theta\theta}\right) - {\tilde\kappa}.
\label{eqn-flowk1}
\end{align}
The details of the discussion in terms of support function can be
referred to \cite{C}.

%% file: study1.tex
\section{Outward Perturbation}
In this section, we will deal with the case that $\varepsilon>0$.
This corresponds to that the initial curve is a small outward
perturbation of an Abresch-Langer self-similar curve.

\begin{lemma}
Let $\varepsilon>0$ and $h(\theta,t)$ be a $\frac{2m\pi}{n}$-periodic
solution to initial value problem~{\em(\ref{eqn-flow0})}.  Then
$$
\max_\theta h(\theta,t) = h(0,t); \qquad
\min_\theta h(\theta,t) = h(m\pi/n,t);
$$
and, $h$ is decreasing for $\theta \in
[0,m\pi/n]$.  Furthermore, let $\kappa$ be the curvature of
the curve supported by $h$, then $\max_\theta \kappa(\theta,t) =
\kappa(0,t)$ and $\min_\theta \kappa(\theta,t)= \kappa(m\pi/n,t)$.
\end{lemma}
\begin{proof}
Firstly, for all $k\geq 1$, $\partial^k_\theta (\hepsilon) =
\left(1+\frac{L}{m\pi}\varepsilon+\varepsilon^2\right)^{-1/2}
\partial^k_\theta(\hmn).$  Therefore, $\hepsilon$ and $\hmn$ have
the same critical points with the same the extremal properties.
As a result,
$$
\max_\theta \hepsilon(\theta) = \hepsilon(0), \qquad\qquad
\min_\theta \hepsilon(\theta) = \hepsilon(m\pi/n)
$$
and $\hepsilon$ decreases for $\theta \in [0,m\pi/n]$.  Secondly,
from the evolution equation,
$$
h_{t\theta} = \frac{1}{(h+\hff)^2}(h_\theta + h_{\theta\theta\theta}).
$$
This shows that $h_\theta$ satisfies a parabolic equation,
$(h_\theta)_t = \kappa^2(h_\theta)_{\theta\theta} + \kappa^2(h_\theta)$.
By Sturm oscillation theorem, \cite{An2}, the number of zeros of $h_\theta$
is non-increasing in $t$.  Thus, for all $t>0$, $h$ has at most 2 critical
points in $[0,m\pi/n]$.

On the other hand, by symmetry of $h$ along $\theta=0$ and $\theta=m\pi/n$,
we must have $h_\theta=0$ at the symmetry.  Hence for each $t$,
$h(\theta,t)$ has exactly one maximum and one minimum for $\theta
\in [0,m\pi/n]$.  The desired results for $h$ follow easily.

The proof for $\kappa$ is similar by simply observing that
$$
\kepsilon = \left( 1 + \frac{L}{m\pi}\varepsilon +
\varepsilon^2 \right)^{1/2} \kmn ;
$$
and the equation for $\kappa_\theta$ is $(\kappa_\theta)_t =
\kappa^2(\kappa_\theta)_{\theta\theta} + 2\kappa\kappa_\theta
(\kappa_\theta)_\theta + 3\kappa^2\kappa_\theta$.
\end{proof}
{\em Remark.\/}
From the above proof, we indeed have some information about
the shape of each leaf of the immersed curve defined by such
$h(\theta,t)$.  In particular,
both the longest distance (the tip of a leaf) and
the shortest distance from the origin to the
curve are attained along the same directions as those of
$\gamma_{m,n}$.

\begin{lemma}\label{lem-maxprinciple}
For sufficiently small $\varepsilon>0$, $\hepsilon(0)<\hmn(0)$
and $\hepsilon(m\pi/n)
> \hmn(m\pi/n)$.  Moreover, $\kepsilon(0)<\kmn(0)$ and
$\kepsilon(m\pi/n) > \kmn(m\pi/n)$.
All these inequalities reverse for $\varepsilon < 0$.
\end{lemma}
\begin{proof}
By definition of $\hepsilon$ and its expansion in $\varepsilon$, we have
$$
\hepsilon = \hmn + \left(1-\frac{L}{2m\pi}\hmn\right)\varepsilon +O(\varepsilon^2).
$$
Note that $\hmn = \kmn$, so
$$
\int_0^{2m\pi} \frac{1}{\hmn} = \int_0^{2m\pi} \frac{1}{\kmn}
= \int_0^{2m\pi}[\hmn + (\hmn)_{\theta\theta}] = \int_0^{2m\pi} \hmn = L.
$$
Applying the Mean Value Theorem, there is a $\theta_*$ such that
$\dfrac{1}{\hmn(\theta_*)}(2m\pi) = L.$  According to the preceding
proposition, $\hmn(0) > \hmn(\theta_*) > \hmn(m\pi/n)$.  As a consequence,
$$
1-\frac{L}{2m\pi}\hmn(0) < 0, \qquad\qquad
1-\frac{L}{2m\pi}\hmn(m\pi/n) > 0.
$$
The inequalities can be obtained by putting the above into the
expansion of $\hepsilon$.

The proof for the inequalities of $\kappa$ is similar by
observing that
$$
\kepsilon = \kmn + \varepsilon\kmn\left(
\frac{L}{2m\pi}-\kmn \right) + \O{\varepsilon^2},
$$
and applying the Mean Value Theorem to $\displaystyle\int_0^{2m\pi}
\kmn = \int_0^{2m\pi} \hmn = L$.
\end{proof}
Thus in the case that $\varepsilon >0$, it follows from the maximum
principle, applied to~(\ref{eqn-flow1}) and~(\ref{eqn-flowk1}), that
there are positive uniform upper and lower bounds for ${\tilde h}$ and
${\tilde\kappa}$ for all $\tau\geq 0$.  As a consequence of parabolic
regularity theory, their higher derivatives are also uniformly
bounded.  From these estimates, we infer the long time existence
of~(\ref{eqn-flow1}) and~(\ref{eqn-flowk1}).

\begin{proposition}\label{prop-hconv}
As $\tau\to\infty$, the solution ${\tilde h}$ to
equation~{\em(\ref{eqn-flow1})} with $\varepsilon>0$ approaches
a stationary solution, which is either $\hmn$ itself or the $m$-circle.
\end{proposition}
\begin{proof}
We consider the entropy $\displaystyle\calE = \int_0^{2m\pi}
\log {\tilde\kappa}$ for the normalized flow~(\ref{eqn-flowk1}).  Then
\begin{align*}
\calE'(\tau) &= \int_0^{2m\pi} \frac{{\tilde\kappa}_\tau}{{\tilde\kappa}}
= \int_0^{2m\pi} \left[{\tilde\kappa} \left({\tilde\kappa}+{\tilde\kappa}_{\theta\theta}\right)-1\right] \\
&= -2m\pi + u(\tau), \\
\intertext{where $\displaystyle u(\tau) = \int_0^{2m\pi} {\tilde\kappa} \left({\tilde\kappa}+{\tilde\kappa}_{\theta\theta}\right).$}
u'(\tau) &= \int_0^{2m\pi} {\tilde\kappa}_\tau \left({\tilde\kappa}+{\tilde\kappa}_{\theta\theta}\right) + {\tilde\kappa} \left({\tilde\kappa}_\tau+{\tilde\kappa}_{\theta\theta\tau}\right) \\
&= 2\int_0^{2m\pi} {\tilde\kappa}_\tau \left({\tilde\kappa}+{\tilde\kappa}_{\theta\theta}\right) \qquad\qquad \text{using integration by parts} \\
&= 2\int_0^{2m\pi} \left[{\tilde\kappa}^2 \left({\tilde\kappa}+{\tilde\kappa}_{\theta\theta}\right) - {\tilde\kappa}\right] \left({\tilde\kappa}+{\tilde\kappa}_{\theta\theta}\right) \\
&= 2\int_0^{2m\pi} {\tilde\kappa}^2 \left({\tilde\kappa}+{\tilde\kappa}_{\theta\theta}\right)^2 - 2\int_0^{2m\pi} {\tilde\kappa} \left({\tilde\kappa}+{\tilde\kappa}_{\theta\theta}\right) \\
&\geq \frac{2}{m\pi} u(\tau) \left(u(\tau)-2m\pi\right).
\end{align*}
Suppose there is a $\tau_1$ such that $u(\tau_1) \geq 2m\pi$.  Then
$u(\tau)$ blows up in finite time and so does $\calE(\tau)$.
This contradicts the boundedness of $\calE$.
Thus, $u(\tau) \leq 2m\pi$ and $\calE$
is decreasing in $\tau$ as in the embedded case, \cite{GH}.

Using this property, we conclude that there is a sequence,
after passing to a subsequence,
$\tau_j\to\infty$ such that $\calE'(\tau_j)\to 0$.  Thus,
using the uniform bounds on ${\tilde\kappa}$, we see that 
${\tilde\kappa}_\tau \to 0$ as $\tau\to\infty$ which in turns
implies ${\tilde\kappa}$ and hence ${\tilde h}$ converge.
\end{proof}

We will rule out that it converges back to an Abresch-Langer curve.
To prepare for this, we first write down some expansions.  Let
$\kepsilon$ be the curvature corresponding to $\hepsilon$.  We have
\begin{align*}
\hepsilon &= \hmn + \varepsilon\left(1-\frac{L}{2m\pi}\hmn\right) + \varepsilon^2 \left(\frac{-1}{2}\hmn - \frac{L}{2m\pi} + \frac{3L^2}{8m^2\pi^2}\hmn\right) + \O{\varepsilon^3}; \\
\kepsilon &= \kmn\left[ 1 + \varepsilon \left(\frac{L}{2m\pi}-\kmn\right) + \varepsilon^2 \left(\frac{1}{2} + \kmn^2 - \frac{L^2}{8m^2\pi^2} - \frac{L}{2m\pi}\kmn \right)\right] + \O{\varepsilon^3}.
\end{align*}
Besides the entropy $\calE$ we have seen before, there is another useful
functional which is also decreasing along the flow.  Let
$$
\calE(\varepsilon) = \int_0^{2m\pi} \log \kepsilon, \qquad\qquad
\calF(\varepsilon) = \int_0^{2m\pi} \log \hepsilon.
$$
Then
\begin{align*}
\calE'(\varepsilon) &= \int_0^{2m\pi} \frac{\kmn}{\kepsilon} \left[ \frac{L}{2m\pi} - \kmn + \varepsilon \left( 1 + 2\kmn^2 - \frac{L^2}{4m^2\pi^2} - \frac{L\kmn}{m\pi} \right) + \O{\varepsilon^2} \right]; \\
\calE''(\varepsilon) &= \int_0^{2m\pi} \frac{\kmn}{\kepsilon} \left( 1 + 2\kmn^2 - \frac{L^2}{4m^2\pi^2} - \frac{L\kmn}{m\pi} \right) - \frac{\kmn^2}{\kepsilon^2} \left( \frac{L}{2m\pi} - \kmn\right)^2 + \O{\varepsilon}. 
\end{align*}
It follows that, using $\kmn=\hmn$,
\begin{align*}
\calE''(\varepsilon) &= 2m\pi - \frac{L^2}{m\pi} + \int_0^{2m\pi} \kmn^2 \\
&= 2m\pi\left[ 1 - 2\left(\meanint \hmn\right)^2 + \meanint \hmn^2 \right]
+ \O{\varepsilon}.
\end{align*}
For $\calF$, we may also work out the same calculations.
\begin{align*}
\calF'(\varepsilon) &= \int_0^{2m\pi} \frac{1}{\hmn} \left[ 1 - \frac{L}{2m\pi}\hmn + \varepsilon \left(-\hmn - \frac{L}{m\pi} + \frac{3L^2}{4m^2\pi^2}\hmn \right) + \O{\varepsilon^2} \right]; \\
\calF''(\varepsilon) &= \int_0^{2m\pi} \frac{1}{\hmn} \left(-\hmn - \frac{L}{m\pi} + \frac{3L^2}{4m^2\pi^2}\hmn \right) - \frac{1}{\hmn^2} \left( 1 - \frac{L}{2m\pi}\hmn \right)^2 + \O{\varepsilon}.
\end{align*}
Using the fact that $1/\hmn = 1/\kmn = \hmn + (\hmn)_{\theta\theta}$,
we have
$$
\calF''(\varepsilon) = -2m\pi \left[ 1 - \left( \meanint \frac{1}{\hmn} \right)^2 + \meanint \frac{1}{\hmn^2} \right] + \O{\varepsilon}.
$$
With these calculations, we have the following nonlinear instability result.
\begin{proposition}\label{prop-hmninstable}
$(\calE + \calF)(\varepsilon) < 0$ for all sufficiently small
$\varepsilon > 0$.
\end{proposition}
\begin{proof}
We consider the functional $\calE + \calF$ and already have
$$
(\calE+\calF)''(\varepsilon) = 2m\pi \left[ \meanint \hmn^2 - \meanint \frac{1}{\hmn^2} \right] + \O{\varepsilon}.
$$
Here,
\begin{align*}
\meanint \frac{1}{\hmn^2} &= \meanint \left[\hmn +(\hmn)_{\theta\theta}\right]^2 \\
&= \meanint \hmn^2 + 2 \meanint \hmn(\hmn)_{\theta\theta} + \meanint (\hmn)_{\theta\theta}^2 \\
&= \meanint \hmn^2 - 2 \meanint (\hmn)_{\theta}^2 + \meanint (\hmn)_{\theta\theta}^2.
\end{align*}
By Poincar\'e Inequality, we have $\displaystyle \meanint
(\hmn)_{\theta\theta}^2 \geq \left(\frac{2\pi}{2m\pi/n}\right)^2
\meanint (\hmn)_{\theta}^2.$  Combining the above, together with
$\dfrac{m}{n} < \dfrac{1}{\sqrt{2}}$, we have
$$
(\calE+\calF)''(\varepsilon) = \left[ 2 - \frac{n^2}{m^2} \right] \meanint (\hmn)_{\theta}^2 + \O{\varepsilon} < 0,
$$
for sufficiently small $\varepsilon$.
\end{proof}
Now, we can finish the proof of part~(a) in the main theorem.  We
have already known that the flow ${\tilde\gamma}$ exists for all
$\tau\geq 0$.  By proposition~\ref{prop-hconv}, it converges either to
the $m$-circle or back to its initial curve.  As $\calE+\calF$ is
non-increasing along the flow, it follows from
proposition~\ref{prop-hmninstable} that the latter is impossible.  So,
it must converge to the $m$-circle and (a)~holds.  Note that we work on
equation~(\ref{eqn-flow1}) which is obtained from normalization by the
algebraic area.  It is equivalent to that by arc length.

%% file: study2.tex
\section{Inward Perturbation}
In this section, we will show that, for $\varepsilon < 0$, the
singularities occur exactly when the area of a leaf vanishes and the
curvature of the curve blows up.  We first establish some analytical
results.  The following lemma basically lays the foundation to show
that in the evolution, the leaves eventually ``shrink'' and exclude
the origin.
\begin{lemma}\label{lem-minmaxh}
Let $\tildeh(\theta,t)$ be a solution to the initial value
problem~{\em(\ref{eqn-flow1})}.  If $\tildeh > 0$ for all $\tau$ in an
interval, then $\tildeh$ is also uniformly bounded above.
\end{lemma}
\begin{proof}
Modifying an argument in \cite{GH},
we consider the following quantity,
$W(\tau)$.
\begin{align*}
W(\tau) &= \int_{-\pi/2}^{\pi/2}
\frac{\cos\theta}{\tkappa(\theta,\tau)} = \int_{-\pi/2}^{\pi/2}
\left(\tildeh + \tildeh_{\theta\theta}\right)\cos\theta \\
&= \int_{-\pi/2}^{\pi/2} \tildeh\cos\theta + \left.
\tildeh_\theta\cos\theta \right|_{-\pi/2}^{\pi/2} + \left.
\tildeh\sin\theta \right|_{-\pi/2}^{\pi/2} - \int_{-\pi/2}^{\pi/2}
\tildeh\cos\theta \\
&= \tildeh(\pi/2) + \tildeh(-\pi/2).
\end{align*}
Therefore, $W(\tau)$ is the width (measured perpendicular to the
longest axis) of a leaf of the curve defined by $h$.
By Jensen's Inequality,
\begin{align*}
\log W(\tau) &\geq \int_{-\pi/2}^{\pi/2} \log\left(
\frac{\cos\theta}{\tkappa(\theta,\tau)} \right) \\
&= \int_{-\pi/2}^{\pi/2} \log\cos\theta - \int_{-\pi/2}^{\pi/2} \log
\tkappa \\
&= \log C - \calE(\tau) \geq \log C - \calE(0),
\end{align*}
for some constant $C>0$, because $\calE$ is decreasing (proved in
proposition~\ref{prop-hconv}).  Thus, $W(\tau)$ is uniformly bounded
below.  Since
$$
W(\tau)\cdot\max_\theta \tildeh(\theta,\tau) \leq \text{the area of a leaf}
\leq A(h) = m\pi,
$$
it follows that $\max_\theta \tildeh(\theta,\tau)$
has a uniform upper bound.
\end{proof}
Now, we may use this lemma to show that for $\varepsilon<0$, at some
time, a curve in the flow will pass through the origin.
\begin{proposition}\label{prop-htouch0}
Let $\tildeh(\theta,\tau)$ be the solution to the initial value
problem~{\em(\ref{eqn-flow1})} with $\varepsilon < 0$.  Then there is a
time $\tau_1$ such that $\min_\theta\tildeh(\theta,\tau_1) = 0$ and
$\tildeh > 0$ for $\tau \in [0,\tau_1)$.
\end{proposition}
\begin{proof}
Suppose on the contrary that $\tildeh > 0$ for all $\tau \geq 0$.
Locally express the
tip of a leaf as a concave graph.  Specifically, we write
$$
\tilde\gamma(\theta(x),\tau) = \left( x, u(x,\tau) \right), \qquad x
\in (-\ell, \ell),
$$
such that $\theta(0) = 0$.  By assumption and the preceding lemma,
there is $M$ such that $\max_x u \leq M$ for all $\tau$.  Note that
the equation on $\tilde\gamma$ corresponding to
equation~(\ref{eqn-flow1}) is
$$
\tilde\gamma_\tau =  \tkappa \bn + \tilde\gamma.
$$
Under this local parametrization, we have
\begin{align*}
\tilde\gamma_\tau &= \pdiff{x}{\tau}\,(1,u_x) + (0,u_\tau); \\
\bn &= \frac{1}{\sqrt{1+u_x^2}}\,(-u_x,1); \\
\tkappa &= \frac{u_{xx}}{(1+u_x^2)^{3/2}}.
\end{align*}
Therefore, we have the parabolic equation
\begin{equation}\label{eqn-flow1grf}
u_\tau = \frac{u_{xx}}{1+u_x^2} - xu_x + u.
\end{equation}
Fixed a small $\xi > 0$, we consider the sub-interval
$(-\ell+\xi,\ell-\xi)$.  By concavity of the graph, we have
$$
\sup\left\{ \modulus{u_x} : x \in (-\ell+\xi,\ell-\xi) \right\} \leq
\frac{1}{\xi} \operatorname{osc}_{(-\ell,\ell)} u \leq \frac{M}{\xi}
$$
Hence, (\ref{eqn-flow1grf}) is uniformly parabolic.  By parabolic
regularity theory, we know that all higher derivatives of $u$ are
uniformly bounded in $(-\ell+\xi,\ell-\xi)$.  In particular, the
curvature at the tip, when $x=0$, is uniformly bounded.  By
propositions~\ref{prop-hconv} and~\ref{prop-hmninstable}, we know that
${\tilde\gamma}$ converges to the $m$-circle, contradicting
lemma~\ref{lem-maxprinciple} and the maximum principle.
\end{proof}
The contradiction shows that $\min_\theta \tildeh(\theta,\tau_1) = 0$.
Clearly, shortly after $\tau_1$, $\min_\theta \tildeh(\theta,\tau) < 0$
and the evolution of the curve enters the third stage in Figure~2b.

Next, we would like to make sure the solution to
problems~(\ref{eqn-flow0}) and~(\ref{eqn-flow1}) exist as long as
the leaves do not shrink to points.  Since the algebraic area may become
negative, it is necessary to look at the unnormalized
equation~(\ref{eqn-flow0}).
\begin{proposition}
Let $h$ be a solution to the initial value problem~{\em(\ref{eqn-flow0})}.
It exists as long as the area enclosed by a leaf does not vanish.
\end{proposition}
\begin{proof}
Suppose otherwise, then there is $t_*>0$ such that actual area of a
leaf~$\geq C > 0$ and $\lim_{t\to t_*^-}\kappa(t) = \infty$.  We
represent the curve at time~$t$, $\gamma(\sbullet,t)$ as a concave
graph as before.  There is
$\delta > 0$ such that for $t \in (t_*-\delta, t_*)$ and
$x\in(-\ell+\xi,\ell-\xi)$, we have the evolution equation
$$
u_t = \frac{u_{xx}}{\sqrt{1+u_x^2}},
$$
which is equivalent to~(\ref{eqn-flow0}); see, for example, \cite{C}.
Again, by similar
argument as above, $u_x$ is uniformly bounded for $x \in
(-\ell+\xi,\ell-\xi)$ and $t\in(t_*-\delta, t_*)$.  Thus, we obtain
a uniform bound for $\kappa(0,t) = \left.
\dfrac{u_{xx}}{\sqrt{1+u_x^2}}\right|_{x=0}$.  By continuity,
as $t\to t_*^-$, $\kappa(0,t) \to$ a finite value, which is a
contradiction.
\end{proof}
\begin{corollary}
Let $t_1$ be corresponding to $\tau_1$ in proposition~{\em\ref{prop-htouch0}}.
There is $t_\infty>t_1$ such that every leaf of the curve shrinks
to a point.
\end{corollary}
\begin{proof}
Since $h(\theta,t_1) > 0$ for all $m\pi/n \ne \theta \in
[0,2m\pi/n]$, we have $A(h(\sbullet,t_1)) > 0$. 
By continuity, there is $\delta
> 0$ such that $A(h(\sbullet,t)) > 0$ for $t < t_1+\delta$.  Therefore,
the area of a leaf of the curve at time~$t < t_1+\delta$ is still
positve.  According to the above proposition, solution
to~(\ref{eqn-flow0}) still exists for $t < t_1+\delta$.  On the other
hand, it is
easy to see that the rate of change of the area of a leaf is less
than~$-\pi$, \cite{AL}.  Thus, there is a time $t_\infty > t_1$ that
all the $m$~leaves disappear.
\end{proof}

{\em Remark.\/}
Recently, an affine version of the curve shortening flow is studied by
B.~Andrews, \cite{Andrews}, also de Lime and Montenegro, \cite{LM}.  In
the latter paper, a classification theorem parallel to the Abresch-Langer
theorem for contracting self-similar curve along the affine flow is
established.  It is interesting to study whether the saddle point
property still holds for these curves.